\documentclass[11pt, a4paper]{article}
\usepackage{amsmath,amssymb,amsfonts}
\usepackage{a4wide}
\usepackage{epsfig,epic,eepic,graphicx}
\newtheorem{theorem}{Theorem}

\newtheorem{prop}[theorem]{Proposition}

\makeatletter

\@addtoreset{equation}{section}
\makeatother

\newcommand{\T}{{\mathbb T}}

\newcommand{\Z}{{\mathbb Z}}
\newcommand{\Q}{{\mathbb Q}}
\newcommand{\R}{{\mathbb R}}

\newcommand{\Sb}{\mathbb{S}}

\newcommand{\pa}{{\partial}}
\newcommand{\na}{{\nabla}}
\newcommand{\dis}{{\displaystyle}}
\newcommand{\eps}{{\varepsilon}}

\def\dis{\displaystyle}

\title{Recent progress in the theory of 
homogenization with oscillating Dirichlet data}
\author{David G\'erard-Varet\footnote{Institut de Math\'ematiques de Jussieu and University Paris 7, 
    175 rue du Chevaleret,75013
   Paris, FRANCE. D. G-V is partially supported by the project  Instabilities in Hydrodynamics} and
Nader Masmoudi\footnote{Courant Institute of Mathematical Sciences, 251 Mercer
   Street, New York, NY 10012, USA.
N. M is  partially supported by
NSF Grant DMS-1211806.} }

\begin{document}
\maketitle

In this talk we study the homogenization of elliptic systems 
 with  Dirichlet boundary
condition, when both the   coefficients  and
the boundary datum are oscillating, namely  $\eps$-periodic. 
In particular, in the paper \cite{GM12},  we showed  that, as $\eps \rightarrow
0$, the solutions  converge  in $L^2$ with a power rate in $\eps$, and
we identified  the homogenized limit system and the homogenized boundary 
data. Due to a boundary 
layer phenomenon, this homogenized system depends in a non trivial way on
the boundary. The 
 analysis in \cite{GM12}   answers a longstanding open  problem,
 raised for instance in \cite{BLP78}.

\section{Introduction} 
 Homogenization of elliptic systems arises in  several 
physical problems where a mixture is present.  
Some of the main applications of the theory are 
 the  diffusion of heat or electricity in a non-homogeneous media,
 the theory of elasticity of mixtures, ...
 Physically, the main goal of the theory   is  to try to compute 
accurate and effective properties of
these mixtures.   Mathematically, we have 
  to find  a limit system towards which 
the   solutions of homogenization problem converge. This 
passage from ``microscopic'' to ``macroscopic'' description 
is called in the literature  ``homogenization''. 


When  both the   coefficients of   the system and
the boundary datum are  oscillating ($\eps$-periodic) and 
due to a boundary 
layer phenomenon, this homogenized system depends in a non trivial way on
the boundary. In this talk, we  answer  a   longstanding open  
problem, raised for instance by Bensoussan, Lions and Papanicolaou 
in their book  ``Asymptotic analysis for periodic structures'' 
  \cite[page xiii]{BLP78}: 
  \begin{quote}
  {\em Of particular importance is the analysis of the behavior of solutions near boundaries and, possibly, any associated boundary layers. Relatively little seems to be known about this problem. }
  \end{quote}

In particular this  result 
 extends   substantially previous works obtained for polygonal
 domains with sides of rational slopes as well as our previous paper 
\cite{GM11} where the case of  irrational slopes was considered. 
We hope that these notes give a better understanding of the proof 
of the result in  \cite{GM12}.

\section{The homogenization problem}

We consider   the homogenization of elliptic systems in divergence form 
\begin{equation}
\label{equation}
  -\na \cdot \left(A\left(\cdot/\eps\right)
\na u\right)(x)    = f,  \quad  x \in \Omega, 
\end{equation}
set in a bounded  domain $\Omega$ of $\R^d$, $\: d \ge 2$, with an
oscillating Dirichlet data 
\begin{equation}
\label{equation2}
u(x)  = \varphi(x,x/\eps), \quad x \in \pa \Omega.
\end{equation}
As is customary,  $\eps > 0$ is a small  parameter,
  and  $A = A(y)$ takes values in $M_d\left(M_N(\R)\right)$, namely 
       $A^{\alpha \beta}(y) \in M_N(\R)$ is a  family of functions of $y \in \R^d$,  indexed by $1\le \alpha,\beta
  \le d$, with values in the set of $N \times N$ matrices. Here,
  $u = u(x)$ and  $\: \varphi = \varphi(x,y)$  take their values in $ \R^N$. We recall, using
  Einstein convention for summation, that for 
each $1\leq i \leq N$, 
$$ (\na \cdot A\left(\cdot/\eps\right)  \na u )_i(x)   \: := \:
  \pa_{x_\alpha} \, \left[
A^{\alpha\beta}_{ij}\left(\cdot/\eps\right)  \, \pa_{x_\beta}  u_j \right](x).  $$
In the sequel,  Greek letters $\alpha, \beta, ...$ will range between 1 and $d$ 
and Latin letters $i,j,k, ...$ will range between 1 and $N$. 

In the context of thermics,  $d=2$ or $3$, $N=1$,  $u$ is the temperature, and  $\sigma = A(\cdot/\eps)
 \na u$ is the heat flux given by Fourier law.  The parameter $\eps$ models heterogeneity, that is
 short-length  variations of the material conducting
 properties. The boundary  term $\varphi$ in
 \eqref{equation2}  corresponds to a prescribed temperature at  the
 surface of the body and $f$ is a source term.
In the context of linear elasticity,  $d=2$ or $3$, $N=d$,  $u$ is the unknown displacement, 
$f$ is the external load and  $A$ is a fourth order tensor that models Hooke's law.

We  make three  hypotheses:
\begin{description}
\item[i)] Ellipticity:
For some $\lambda > 0, \:$  for  all family of vectors $\xi = \xi^\alpha_i
\in \R^{Nd}$ 
\begin{equation*} \label{ellipticity}
 \lambda \,\sum_{\alpha}  \xi^\alpha \cdot \xi^{\alpha} \: \le \:
 \sum_{\alpha, \beta, i, j } A^{\alpha,\beta}_{ij}
 \, \xi^\beta_j \,  \xi^\alpha_i \: \le \: \lambda^{-1}  \sum_{\alpha} \xi^\alpha
\cdot \xi^{\alpha}.
\end{equation*}
\item[ii)] Periodicity: $\displaystyle  \: \forall y \in \R^d, \: \forall h \in \Z^d, \: \forall x \in \pa \Omega, \:\: A(y+h)
  \: = \: A(y),  \quad \varphi(x,y) \: = \:  \varphi(x,y+h)$. 
\item[iii)] Smoothness: The functions $A, f$ and  $\varphi$, as well as the domain
  $\Omega$ are smooth. It is actually enough to assume that $\varphi$ and $\Omega$ are in some 
 $H^s$ for $s$ big enough, but we will not try to compute the optimal regularity.    
\end{description}

The main question we are trying to answer is the following: 

\medskip
\underline{Question}: What is the limit
 behavior of the solutions $u^\eps$  as $\eps \rightarrow 0$ ?  
Can we go beyond the limit and compute a full expansion of $u^\eps$ ?
 
\medskip 
This question goes back at least to the 1970's, and a classical approach consists in trying a two-scale expansion: 

\medskip
\underline{Classical approach}: {\em Two-scale asymptotic expansion}: 

\begin{equation} \label{exp}
\boxed{ u^{\eps}_{app} \: = \:  u^0(x) \: + \: \eps u^1(x,x/\eps)  \: + \:  \dots \: + \: \eps^n u^n(x,x/\eps)} 
\end{equation}
with  $u^i = u^i(x,y)$ periodic in $y$.

\section{Case without boundary}

The two-scale approach works well in the case without boundary, namely in the 
whole space case or in the case of a periodic domain (say of period 
1 and $\eps$ is taken to be equal to  $1/n$ with $n$ an integer). 
In particular one can construct  inductively all the terms in 
the expansions. Let us recall few classical facts  (see for instance 
\cite{Sanchez80,MT97,JKO94,CD99b}) :

\begin{description}

\item[i)]  The construction of the $u^i$'s  involves the famous {\em cell problem}  
\begin{equation} \label{chi-eq}
\boxed{ -  \na \cdot  \: \left( A \na \chi^\gamma \right)(y)   \: = \: 
 \na_{\alpha} \cdot A^{\alpha \gamma} \, (y), \quad y  \mbox{ in } \T^d} 
\end{equation}
 with solution  $\chi^\gamma  \in M_N(\R) $.

\item[ii)]  The solvability condition for $u^2$ yields
 the equation satisfied by $u^0$, namely 
 $u^0$ (which does not depend on $y$) 
 satisfies
\begin{equation} \label{uf}
\na \cdot A^0 \na  u^0  \: = \: f  
\end{equation}
where the constant homogenized matrix is given by 
$$ \boxed{ A^{0,\alpha\beta} \: = \: \int_{\T^d} A^{\alpha\beta} (y) \, dy  \: +
\: \int_{\T^d} A^{\alpha\gamma} (y) \partial_{y_\gamma} \chi^\beta(y) \, dy. }$$

\end{description}

The second term in the expansion \eqref{exp} reads
\begin{equation} \label{u1}
 u^1(x,y) \: := \: \tilde{u}^1(x,y) +  \bar{u}^1(x)  \:  := \:
-\chi^\alpha(y)  \pa_{x_\alpha} u^0(x) + \bar{u}^1(x),
\end{equation}
where $\chi$ is again the solution of \eqref{chi-eq}.

To find an equation for the  average part $\bar{u}^1(x)$,   
 one needs to
 introduce another family
 of 1-periodic matrices
 $$\Upsilon^{\alpha\beta} = \Upsilon^{\alpha\beta}(y)  \in M_n(\R), \: \alpha,\beta
= 1, ...,d,$$
 satisfying
\begin{equation} \label{upsilon}
 -\na_y \cdot  A \, \na_y \Upsilon^{\alpha\beta}\:  = \: B^{\alpha\beta} - \int_y
B^{\alpha\beta}, \quad \int_y \Upsilon^{\alpha\beta} = 0,
\end{equation}
where
$$ B^{\alpha\beta}\: := \: A^{\alpha\beta} - A^{\alpha\gamma}
\frac{\pa \chi^\beta}{\pa y_\gamma} - \frac{\pa}{\pa y_\gamma} \left(
  A^{\gamma \alpha} \chi^\beta \right). $$
Formal considerations yield
\begin{equation} \label{u2}
u^2(x,y) \: := \: \Upsilon^{\alpha,\beta} \frac{\pa^2 u^0}{\pa x_\alpha \pa
  x_\beta} \: - \: \chi^\alpha \pa_\alpha \bar{u}^1 +  \bar{u}^2  
\end{equation}
and that the average term  $\bar{u}^1=\bar{u}^1(x)$  formally satisfies
the equation
\begin{equation} \label{averageu1}
 -\na \cdot A^0 \na \bar{u}^1 = c^{\alpha\beta\gamma} \frac{\pa^3
   u^0}{\pa x_\alpha \pa x_\beta \pa x_\gamma}, \quad c^{\alpha\beta\gamma}
 \: := \: \int_y A^{\gamma\eta} \frac{\pa
  \Upsilon^{\alpha\beta}}{\pa y_\eta} - A^{\alpha\beta} \chi^{\gamma}.
\end{equation}
We refer to \cite{AA99} for more  details.

Inductively, one can keep constructing all the terms 
of the expansion by introducing new corrector families 
as in \eqref{upsilon} and solving homogenized systems 
to determine $\bar{u}^k  $ as in \eqref{averageu1}. 
Note that in this case, we do not need an extra boundary condition 
to solve \eqref{averageu1}.

\section{Case with   boundary}

 Two boundary conditions   have been widely studied and
are by  now well understood as long as we are only interested in 
the first term of the expansion:
\begin{enumerate}
\item {\em The non-oscillating Dirichlet problem}, that is
  \eqref{equation} and \eqref{equation2} with $\varphi = \varphi(x)$.
\item {\em The oscillating Neumann problem}, that is
  \eqref{equation} and 
\begin{equation} \label{equation3}
n(x) \cdot \left( A(\cdot/\eps) \na u \right)(x) \: = \: \varphi(x,x/\eps), \quad x
\in \pa \Omega, 
\end{equation} 
where $n(x)$ is the normal vector
and 
 with a standard
  compatibility condition on  $\varphi$. Note that 
 in thermics, this boundary condition  corresponds to a
given heat flux  at the solid surface.
\end{enumerate} 
Notice  that in both  problems, the usual energy estimate provides a uniform bound on
the solution $u^\eps$ in $H^1(\Omega)$. 
\medskip

For the non-oscillating Dirichlet problem, one shows that
$u^\eps$ weakly converges in $H^1(\Omega)$ to the solution $u^0$ of the homogenized
system 
\begin{equation} \label{homogene}
\left\{
\begin{aligned}
  -\na \cdot \left( A^0 \na u^0\right)(x)  = f, & \quad  x \in \Omega,\\
 u^0(x)  = \varphi(x), & \quad x \in \pa \Omega.
\end{aligned}
\right.
\end{equation}
It is also  proved in \cite{BLP78} that 
\begin{equation} \label{ansatz}
u^\eps(x) \: = \:  u^0(x) \: + \: \eps u^1(x,x/\eps) \: + \:
O(\sqrt{\eps}), \: \mbox{ in  } \: H^1(\Omega).  
\end{equation}
Actually, an open problem in this area was to compute the next term in 
the expansion in the presence of a boundary, namely to compute 
 $ u^1(x,x/\eps)  $. Indeed, it is not difficult to see that 

\begin{equation} \label{u1*}
 u^1(x,y) \: = \:
-\chi^\alpha(y)  \pa_{x_\alpha} u^0(x) + \bar{u}^1(x),
\end{equation}
where  $\bar{u}^1(x) $ solves the homogenized equation 
\eqref{averageu1}. However, the main difficulty is to find 
the boundary data for $\bar{u}^1(x) $.
The new analysis of \cite{GM12} gives an answer to this problem
 (see also next section). 

\medskip
For the oscillating Neumann problem, two cases must be distinguished.
On one hand, if  $\pa \Omega$ does not contain flat pieces, or if it
contains finitely many flat pieces 
whose normal vectors do not belong to $\R  \Z^n$, then 
$$   \varphi(\cdot, \cdot/\eps) \rightarrow \overline{\varphi} := \int_{[0,1]^d} \varphi
\: \mbox{ weakly in} \:   L^2(\pa \Omega)$$ 
and $u^\eps$ converges weakly to the solution $u^0$ of 
\begin{equation} \label{homogene2}
\left\{
\begin{aligned}
  -\na \cdot \left( A^0 \na u^0\right)(x)  = 0, & \quad  x \in \Omega,\\
  n(x) \cdot \left( A^0 \na u^0 \right)(x) =  \overline{\varphi}(x), & \quad x \in \pa \Omega.
\end{aligned}
\right.
\end{equation}
On the other hand, if  $\pa \Omega$ does contain a flat piece whose normal vector
  belongs to  $\R \, \Q^d$,
  then the family  $\varphi(\cdot, \cdot/\eps)$ may have
  a continuum of accumulation points as $\eps \rightarrow 0$. Hence,  $u^\eps$
  may have a continuum of accumulation points in $H^1$ weak, corresponding
  to  different Neumann boundary data. We refer to \cite{BLP78} for all
  details.

\section{Case of an oscillating    Dirichlet data} \label{bd}

Here we study \eqref{equation} with the boundary data  \eqref{equation2}.  
One of the motivation to study this case is actually to 
understand the  boundary condition for $\bar{u}^1(x) $ which 
appears in  \eqref{u1}. \medskip

Let us explain the two main sources of difficulties in  
studying \eqref{equation}-\eqref{equation2}:

\begin{description}
\item[i)]
 One has  uniform $L^p$ bounds on the solutions $u^\eps$ of
 \eqref{equation}-\eqref{equation2}, but  no  uniform $H^1$ bound {\it a
   priori}. This is due to the fact that  
$$\| x \mapsto \varphi(x,x/\eps) \|_{H^{1/2}(\pa \Omega)} = O(\eps^{-1/2}), 
 \: \mbox{ resp.} \:   \| x \mapsto \varphi(x,x/\eps) \|_{L^p(\pa
   \Omega)} = O(1), \: p >1. $$    
 The usual energy  inequality, resp. the estimates in article
\cite[page 8, Thm 3]{AvLi} yields   
$$ \| u^\eps \|_{H^1(\Omega)} =  O(\eps^{-1/2}),  \: \mbox{ resp.} \:  \| u^\eps \|_{L^p(\Omega)} =  O(1),
 \: p > 1.$$

This indicates that singularities of $u^\eps$ are {\it a priori} stronger than in the
usual situations. It is rigorously established in the core of the
paper \cite{GM12}. 
\item[ii)] Furthermore, one can not expect these stronger singularities
  to be periodic oscillations. Indeed, the oscillations of $\varphi$ are
   at the boundary, along which they do not have any periodicity
  property. Hence, it is reasonable that $u^\eps$ should exhibit
  concentration  near $\pa \Omega$, with no periodic character, as $\eps
  \rightarrow 0$. This is a so-called {\em boundary layer
    phenomenon}. The key point is to describe this boundary layer,
  and  its effect on the possible weak limits of $u^\eps$.

 \end{description}

It is important to note that  there is also a  boundary layer in the 
non-oscillating Dirichlet problem, although it has in this case a
lower amplitude (it is only necessary to compute the 
boundary data of $\bar{u}^1$ to solve 
\eqref{averageu1}). More precisely, it is responsible for the 
$O(\sqrt{\eps})$ loss in the error estimate \eqref{ansatz}. If 
either the $L^2$ norm,  or the $H^1$ norm in  a relatively compact   subset $\omega \Subset \Omega$ is considered, one may avoid this loss  as strong gradients near the boundary are filtered out. Following  Allaire and Amar (see \cite[Theorem 2.3]{AA99}), we can give a more precise description than \eqref{ansatz}:
\begin{equation} \label{ansatz2}
u^\eps  \: = \:   u^0(x) \: + \: O(\eps)  \: 
\mbox{ in } \:  L^2(\Omega), \quad   u^\eps(x)  \: =  \:  u^0(x) + \eps u^1(x,x/\eps) \: + \: O(\eps)  \: \mbox{ in } \: H^1(\omega).
\end{equation} 
Still following \cite{AA99}, another way to put the emphasis on  the
boundary layer is to introduce  the solution $u^{1,\eps}_{bl}(x)$ of 
\begin{equation} \label{u1bl}
\left\{
\begin{aligned}
  -\na \cdot A\left(\frac{x}{\eps}\right) 
\na u^{1,\eps}_{bl}  = 0, & \quad  x \in \Omega \subset \R^d,\\
 u^{1,\eps}_{bl}  = -u^1(x,x/\eps), & \quad x \in \pa \Omega.
\end{aligned}
\right.
\end{equation} 
Actually, understanding  this system and requiring that  $u^{1,\eps}_{bl} $  goes to 
zero inside the domain $\Omega$ allows to determine the right boundary condition 
for $\bar{u}^1 $. 
Hence,  one can show that  
\begin{equation} \label{ansatz3}
u^\eps(x)  \: =  \:  u^0(x) + \eps u^1(x,x/\eps) + \eps
u^{1,\eps}_{bl}(x) \: + \: O(\eps),  \: 
\mbox{ in } \:  H^1(\Omega).
\end{equation}
or
\begin{equation} \label{ansatz4}
u^\eps(x)  \: =  \:  u^0(x) + \eps u^1(x,x/\eps) + \eps
u^{1,\eps}_{bl}(x) \: + \: O(\eps^2),  \: 
\mbox{ in } \:  L^2(\Omega).
\end{equation}

Note that system \eqref{u1bl} is a special case of
\eqref{equation}-\eqref{equation2}. Thus, the homogenization  of the
oscillating Dirichlet problem may give a refined description of the
non-oscillating one.

\section{Prior results} \label{pr}

 Until recently, {\em results  were all limited to convex polygons
with rational normals.} This  means that  
\begin{equation*}
\Omega \: := \:  \cap_{k=1}^K   \left\{ x, \quad n^k \cdot x  >   c^k \right\}
\end{equation*}
is bounded by $K$ hyperplanes, {\em  whose unit normal vectors $n^k $  
 belong to  $\R \, \Q^d$.} Under this  assumption, the study  of
\eqref{equation}-\eqref{equation2} can be carried out. The
keypoint is  the addition of boundary layer correctors to the formal
two-scale expansion:
\begin{equation}
 u^{\eps}(x) \: \sim  \: u^0(x) \: + \:  \eps u^1(x,x/\eps) \: + \: \sum_k
v_{bl}^k\left(x,\frac{x}{\eps}\right),  
\end{equation}
where $v_{bl}^k = v_{bl}^k(x,y) \in \R^n$ is defined for $x \in \Omega$, and $y$
in the half-space
$$\Omega^{\eps,k} \: = \:
\left\{ y , \quad   n^k \cdot y > c^k/\eps    \right\}.$$
These correctors satisfy 
\begin{equation} \label{BL1a}
\left\{
\begin{aligned}
 -\na_y \cdot \, A(y) \na_y \, v_{bl}^k  = 0, & \quad y \in \Omega^{\eps,k},\\
   v_{bl}^k  = \varphi(x,y) - u^0(x), & \quad  y \in \pa \Omega^{\eps,k}.
\end{aligned}
\right.
\end{equation} 
We refer to the papers by Moskow and Vogelius \cite{MV97}, and Allaire
and Amar \cite{AA99} for more details. These papers deal with the special case
\eqref{u1bl}, but the results adapt to more  general
oscillating data. Note that $x$ is just  a parameter in
\eqref{BL1a} and that the assumption $n^k \in
\R \, \Z^d$ yields periodicity of the function $A(y)$ tangentially to the
hyperplanes. The periodicity property is used in a crucial way in the aforementioned references. First, it yields easily well-posedness of the
boundary layer systems \eqref{BL1a}. Second, as was  shown by Tartar in \cite[Lemma 10.1]{JLL} (see also subsection \ref{rat}), the solution  $v_{bl}^k(x,y)$ converges exponentially fast  to some $v_{bl,*}^k(x)= \varphi_*^k(x) - u^0(x)$, when $y$ goes to infinity
transversely to the $k$-th hyperplane. In order for the boundary
layer correctors to vanish at infinity (and to  be $o(1)$ in $L^2$),
one must  have $v_{bl,*}^k = 0$, which  provides the  boundary condition
for $u^0$. Hence, $u^0$ should satisfy a system of the type 
\begin{equation} \label{homogene3}
\left\{
\begin{aligned}
  -\na \cdot \left( A^0 \na u^0\right)(x)  = f, & \quad  x \in \Omega, \\
 u^0(x)  = \varphi_*(x), &\quad x \in \pa \Omega. 
\end{aligned}
\right.
\end{equation}
where $\varphi_*(x) \: := \: \varphi_*^k(x)$ on the $k$-th
side of $\Omega$. Nevertheless,  this  picture is not completely correct. Indeed, 
 there is
still {\it a priori} a dependence  
of $\varphi_*^k$ on $\eps$,  through the domain
$\Omega^{\eps,k}$. In fact, Moskow and Vogelius exhibit 
 examples for which there is an infinity of accumulation points for
 the $\varphi_*^k$'s, as $\eps \rightarrow 0$. Eventually, they show that the
 accumulation points of $u^\eps$ in $L^2$ 
 are the solutions $u^0$ of systems  like \eqref{homogene3}, in which the
 $\varphi_*^k$'s are replaced by
 their  accumulation points. See \cite{MV97} for rigorous
 statements and proofs. We stress that their  analysis
 relies heavily  on the special  shape of $\Omega$, especially the rationality assumption. 

\medskip
A step towards more generality has been made in our recent paper
\cite{GM11} (see also \cite{GM10}),
 in which generic convex polygonal domains are
considered. Indeed, we assume in \cite{GM11} that {\em the normals $n = n^k$ satisfy the
Diophantine condition:}
\begin{equation} \label{diophantine}  
 \mbox{ For all } \:   \xi \in \Z^d  \setminus \{0\}   \quad  |P_{n^\bot}(\xi)|  >
\kappa \, |\xi|^{-l}, \quad \mbox{for some } \:  \kappa, \: l > 0,  
\end{equation}
 where $P_{n^\bot}$ is the projector orthogonal to 
$n$.  Note  that for dimension $d=2$ this condition  amounts to:
\begin{equation*} 
 \mbox{ For all } \:   \xi \in \Z^d  \setminus \{0\}   \quad 
   | n^\bot \cdot \xi | \:  := \:   |-n_2 \xi_1 + n_1 \xi_2| > \kappa \, |\xi|^{-l}, \quad \mbox{for some } \:  \kappa, \: l > 0,  
\end{equation*}
whereas for $d=3$, it is equivalent to:
\begin{equation*} 
 \mbox{ For all } \:   \xi \in \Z^d  \setminus \{0\}   \quad 
 |\dis  n \times  \xi| > \kappa \, |\xi|^{-l}, \quad \mbox{for some } \:  \kappa, \: l > 0.
\end{equation*}
Condition \eqref{diophantine} is generic in the sense that it holds for  almost every 
$n \in S^{d-1}$. 

\medskip
Under this Diophantine assumption,  one can perform the homogenization of problem  \eqref{equation}-\eqref{equation2}. {\em Stricto sensu},  only the case \eqref{u1bl}, $d=2,3$ is treated in \cite{GM11},  but our analysis extends straightforwardly to the general setting.  
Despite a loss of periodicity in the tangential variable, we manage to
solve the boundary layer equations, and prove convergence of $v^{k}_{bl}$
 away from the boundary. The main idea is to work with quasi-periodic functions instead 
of periodic ones (see also subsection \ref{nrat}). 
 Interestingly, and contrary to the
``rational case'', the field $\varphi_*^k$ does not depend on $\eps$. As a result,
we establish convergence of the whole sequence $u^\eps$ to the single
solution $u^0$ of \eqref{homogene3}. We stress that, even in this polygonal
setting, the boundary datum $\varphi_*$ depends in a non trivial way on
the boundary. 
 In particular, it is not simply the average of $\varphi$ with
respect to $y$, contrary to what happens in the Neumann case.

\section{Main new  result and sketch of proof}


The main new result of \cite{GM12} is to treat the case 
of a   smooth domain: 
\begin{theorem}   {\bf (Homogenization in smooth domains) }  \label{theo}

\smallskip
\noindent
Let $\Omega$ be a smooth bounded domain of $\R^d$, $d \ge 2$.  We assume that it is 
uniformly  convex (all the principal curvatures are bounded from below).    

\smallskip
\noindent
Let $u^\eps$ be the solution of system \eqref{equation}-\eqref{equation2}, under the ellipticity, periodicity
and smoothness conditions i)-iii). 

\smallskip
\noindent
There  exists a boundary term $\varphi_*$ (depending on $\varphi$, $A$ and $\Omega$), with  $\varphi_*  \in L^p(\pa \Omega)$  for all finite $p$, and a solution 
$u^0$ of \eqref{homogene3}, with  $u^0 \in L^p(\Omega)$ for all finite $p$, such that:   
\begin{equation} \label{rate}
\| u^\eps - u^0 \|_{L^2(\Omega)} \: \leq  \: C_\alpha \, \eps^\alpha, \quad \mbox{ for all } \: 0 < \alpha < \frac{d-1}{3d+5}. 
\end{equation}
\end{theorem}

We will present a sketch of the proof of theorem \ref{theo}: 

>From the two difficulties explain in section \ref{bd}, we know that 
 the first term in the expansion \eqref{exp} should 
be independent of $y$ and should solve \eqref{uf}. The main question is :

\underline{Question}: What is the boundary value $\varphi^0$ of  $u^0$ ? 

\underline{Solution}:   We need a  boundary layer corrector

\underline{Difficulty}:  There is no clear structure for the boundary layer. 

\underline{Guess}: The boundary layer has typical scale $\eps$ and 
there are  no curvature effect: 

\begin{itemize} 
\item Near a point $x_0 \in \pa \Omega$, we  replace $\pa \Omega$ by the tangent plane at $x_0$:
$$ T_0(\pa \Omega) \: := \: \{ x, \: x \cdot n_0 = x_0 \cdot n_0  \}  $$ 
\item  We dilate by a factor $\eps^{-1}$.
\end{itemize}

\medskip
Formally, for $x \approx x_0$, one looks for 
$$\boxed{ u^{\eps,bl}(x) \: \approx \: U_0(x/\eps) }$$ 

where the profile  $U_0 = U_0(y)$ is defined in the half plane 
$$ H^{\eps}_0 \: = \: \{y,  \: y \cdot n_0 >  \eps^{-1} x_0 \cdot n_0 \}.$$

\medskip
It satisfies the system: 
 \begin{equation}  \label{U0} 
\boxed{\left\{
\begin{aligned}
 \na_y \cdot (A \na_y U_0 )  = 0 & \quad \mbox{in } H^{\eps}_0, \\  
 U_0\vert_{\pa H^\eps_0}   = \varphi - \varphi^0(x_0).&  
\end{aligned}
\right.}
\end{equation}

Notice that in this system, $x_0$ is just a parameter.

\subsection{Study of an auxiliary boundary layer system }
The previous heuristic justifies the study of 
\begin{equation} \label{BL} \tag{BL} 
\boxed{\left\{
\begin{aligned}
 \na_y \cdot (A \na_y U )  = 0 & \quad \mbox{in } H, \\  
 U\vert_{\pa H}     = \phi. &  
\end{aligned}
\right.}
\end{equation}
where $H \: := \: \{ y, \quad y \cdot n >  a \} $  and 
$\phi$ is  $1$-periodic in $y$. 

We expect that the solution  $U$ of \eqref{BL} satisfies:  
$$ U \rightarrow U_{\infty} (\phi), \quad \mbox{ as } y \cdot n  \rightarrow +\infty,$$ 
 for some constant  $U_\infty =  U_{\infty} (\phi) $ that depends linearly on $\phi$. 

If we go back to  $U_0$ which solves \eqref{U0}, 
one can derive the  homogenized boundary data $\varphi^0.$ Indeed: 
\begin{itemize}

\item On one hand, one wants $U_0 \rightarrow 0$ (localization property) 
when   $  y \cdot n  \rightarrow +\infty $.

\item On the other hand, 
$$ U_0 \rightarrow U_\infty (\varphi - \varphi^0(x_0)) \: = \: U_\infty(\varphi) - \varphi^0(x_0) $$
so that we need to take: 
$$ \varphi^0(x_0) \: := \: U_\infty(\varphi) .   $$
\end{itemize} 
    

 This formal reasoning raises many problems : 

\begin{enumerate}
\item {\em The well-posedness of \eqref{BL} is unclear: }

\medskip
- No natural functional setting (no decay  along the boundary). \\
- No Poincar\'e inequality. \\
- No maximum principle. 

\medskip
\item {\em  The existence of a limit $U_\infty$ for \eqref{BL}  is unclear: } 
 
There is an underlying problem of ergodicity. 

\medskip
\item  {\em$U_\infty$ depends also on $H$, that is on $n$ and $a$:} 
 
\medskip
- There is  no obvious  regularity of $U_\infty$ with respect to $n$.

\medskip
- Back to the original problem,  our definition of  $\varphi^0(x_0)$  depends on $x_0$, but also on the subsequence  $\eps$. 
Indeed, there is possibly many accumulation points as $\eps \rightarrow 0$ (see \cite{MV97}). 

\end{enumerate}

\subsection{Polygons with sides of rational slopes} \label{rat}

In this  cases, the  boundary layer systems of 
type \eqref{BL} can be fully understood  (see \cite{MV97,AA99}). 
For simplicity, we only concentrate on the case $d=2$. 

\medskip
\begin{enumerate}
\item Well-posedness: {\em The coefficients of the systems are periodic tangentially to the boundary}. 
After rotation, they turn into systems of the  type 
\begin{equation} \label{BL1} \tag{BL1} 
\boxed{\left\{
\begin{aligned}
 \na_{z} \cdot (B \na_{z} V )  = 0, & \quad z_2 > a,  \\  
 V\vert_{z_2 = a }   = \psi, &  
\end{aligned}
\right.}
\end{equation}
with coefficients and boundary data that are  periodic in  $z_1$ which 
yields a natural variational formulation. 

\item Existence of the limit : {\em   Saint-Venant estimates} on \eqref{BL1}. 

One shows that  $ F(t) \: := \: \int_{z_2 > t} | \na_z V |^2 \, dz  $
satisfies the differential inequality.
$$  { F(t) \: \le \: - C F'(t)}.  $$ 
>From there, one gets exponential decay of all derivatives, and the 
fact that: 
$$ V \rightarrow V_\infty, \: \mbox{ exponentially fast, as } \: z_2 \rightarrow +\infty $$
 and hence   going back to $(BL)$, we get
$$ U  \rightarrow U_\infty, \: \mbox{ exponentially fast, as } \: y \cdot n  \rightarrow +\infty.  $$

A {Key} ingredient in this case is the 
 Poincar\'e  inequality  for functions periodic in $z_1$ with zero mean.

\item In polygonal domains, the regularity of $U_\infty$ with respect to $n$ does not matter.  However, for 
rational slopes, the limit $U^\infty$ does depend on $a$.  
This means that if we go back to our original problem 
 (in polygons with rational slopes),  The analogue of our theorem is only available up to subsequences in $\eps$. 
 Moreover,  the boundary data of the 
 homogenized system may depend  on the subsequence.
Indeed, there are  examples with a continuum of accumulation points 
(see \cite{MV97}). 
\end{enumerate}

\subsection{More general treatment of \eqref{BL} } \label{nrat}

It is worth pointing out that one 
can not be fully general: 
The existence of $U_\infty$ requires some ergodicity property. 
A simple example is : 
 $$\boxed{\: \Delta U = 0 \quad \mbox{in } \{ y_2
  > 0 \}, \quad U \vert_{y_2=0} = \phi}.$$  

\medskip
\begin{itemize}
\item If $\phi$ 1-periodic, then $U(0,y_2) \rightarrow  \int_0^1 \phi$
  exponentially fast.    

\medskip
\item {\em  But there exists  $\phi \in L^\infty$ such that 
  $U(0,y_2)$ has  no limit.}  
\end{itemize}

Indeed, we have an explicit formula:  $\displaystyle \:  U(0,y_2) \: = \: \frac{1}{\pi} \, \int_\R \frac{y_2}{y_2^2+ t^2} \phi(t) \, dt. $
  For $\phi$ with values in $\{ +1, -1\}$, the asymptotics relates to  
{\em  coin
  tossing}.   Hence, we need some extra structure (or ergodicity) 
to solve the problem. 

In our case, we have some ergodicity property !  For general half planes, the coefficients of \eqref{BL} or \eqref{BL1} are not periodic, {\em  but they are  quasiperiodic in  the tangential variable.} 
We recall that a  function $F = F(z_1)$ is quasiperiodic if it reads
$$\boxed{ F(z_1) = {\cal F}(\lambda z_1)}, $$ 
 where $\lambda \in \R^D$ and  ${\cal F} = {\cal F}(\theta)$ is  periodic over $\R^D$ ($\: D \ge 1$).   
As an {example}: For \eqref{BL1}, $\: D = 2 \: $,  and $\: \lambda = n^\bot$ (the  tangent vector). 

Notice that the {previous results} (subsection \ref{rat})
 correspond to the case: 
$n \in \R \, \Q^2$. 
Now, we  replace this by the small divisor assumption: 

\begin{center}
\fbox{ (H) {\em  $\quad \exists \kappa > 0$,   $\displaystyle  | n \cdot \xi | \ge \kappa |\xi|^{-2}, \quad \forall \xi \in \Z^2\setminus\{0\}. $} }
\end{center}

Note that the  assumption (H) is generic in the normal $n$: It is 
 satisfied by  a set of full measure in $\mathbb{S}^1$.  
But it does not include the previous result of subsection \ref{rat}.

  
\begin{prop}
  {\em If $n$ satisfies (H), the system \eqref{BL} is "well-posed", with a smooth solution $U$ that converges fast to some constant  $U_\infty$. 
  Moreover, $U_\infty$ does not depend on $a$. }
  \end{prop} 
 
\underline{Proof of the proposition}:
 \begin{enumerate}
 \item Well-posedness: involves quasiperiodicity. One has: 
\begin{equation*} 
\boxed{\left\{
\begin{aligned}
 \na_{z} \cdot (B \na_{z} V )  = 0, & \quad z_2 > a,  \\  
 V\vert_{z_2 = a }   = \psi, &  
\end{aligned}
\right.}
\end{equation*}
 where  $\displaystyle \: B(z) \: = {\cal B}(\lambda z_1, z_2), \:  \psi(z) \: = \: {\cal P}(\lambda z_1, z_2)$.  

Notice that the functions  ${\cal B} = {\cal B}(\theta, t)$ and  ${\cal P} = {\cal P}(\theta,t)$ are periodic in $\theta \in \T^2$. 
   
The idea is to  consider an enlarged system in $\theta, t$, of unknown ${\cal V} = {\cal V}(\theta,t )$: 
\begin{equation} \label{BL2} \tag{BL2} 
\boxed{\left\{
\begin{aligned}
D  \cdot ({\cal B} D {\cal V})  = 0, & \quad  t> a,  \\  
 {\cal V}\vert_{t = a} = {\cal P}  &  
\end{aligned}
\right.}
\end{equation}
where $D$ is the "degenerate gradient" given by $ \: D = (\lambda \cdot \na_\theta, \pa_t) $

\medskip
\underline{Advantage}: Back to a periodic setting ($\theta \in \T^2$). 

\medskip 
\underline{Drawback}: We have a  degenerate elliptic equation.  However, 
we are still able to prove the following : 

\medskip
- Variational formulation with  a unique weak solution ${\cal V}$. 

\medskip
- One can prove through energy estimates than ${\cal V}$ is smooth. 

\medskip
- Allows to recover $V$ through the formula $\: V(z) \: = \: {\cal V}(\lambda z_1, z_2)$.

\item To prove the  convergence to a constant at infinity, we 
rely again  on Saint-Venant type  estimates, adapted  to \eqref{BL2}. 
Thanks to (H), we prove that $\displaystyle F(t) \: := \int_{t' > t} | D {\cal V}|^2 \, d \theta \, dt' $ satisfies 
$$ \boxed{F(t) \: \le \: C(-F'(t))^{\alpha}, \quad \forall \alpha < 1. } $$
But, we have only polynomial convergence towards a constant. 
\end{enumerate}

We point out that this better understanding of the auxiliary boundary layer systems allows to handle the  generic polygonal domains in the next 
subsection.

\subsection{ Extension to smooth domains}

{The are at least three  main difficulties 
 to extend the previous analysis to smooth domain} : 
 
\begin{enumerate}
\item    The  none  smoothness of $U_\infty$ with respect to $n$.  
Indeed,   $U_\infty$    is only defined almost everywhere (diophantine assumption).   
  
\underline{Idea}:  {\em For any $\kappa > 0$, we can prove that 
 $U_\infty$ is Lipschitz when it is  restricted  to 
 $$ A_\kappa \: := \: \left\{ n \in \Sb^1, \: | n \cdot \xi | \: \ge \: \frac{\kappa}{|\xi|^2}, \: \forall \xi \in \Z^2\setminus\{0\} \right\}.$$ }
Moreover, we have that  $\left| A_\kappa^c \right| = O(\kappa)$. 

In the course of the proof,  the construction of the boundary layer corrector can be performed in 
 the vicinity of points $x$ such that $n(x) \in A_\kappa$.  In some sense, the 
contribution of the remaining part of the boundary is negligible when $\kappa \ll 1$.  More precisely, 

\item  We have to approximate the  smooth domains  by some polygons with sides having 
normal vectors in the set  $ A_\kappa$. In doing so, we will introduce another small 
parameter  $\eps^\alpha$. 

\item  We have to construct a more accurate approximation due to the many errors made 
 in the previous two points.

\end{enumerate}

Broadly, optimizing in $\kappa$, $\alpha$ and $\eps$  yields a rate of convergence. 
We refer to \cite{GM12} for the details.

\section{Conclusions}

We would like to conclude by mentioning  a few related results. 
Recently there was many activity in the theory of homogenization
and many new problems were addressed. We would like to mention some of 
them since we think they may give a better understand of our result 
or/and   may  be combined with our result: 
\begin{itemize} 
\item Our results on the boundary data problem were recently extended to the eigenvalue problem, see \cite{Prange1}. Also, the behavior of the reduced boundary layer system \eqref{BL} was recently investigated by C. Prange in \cite{Prange2}, without any diophantine assumption. 

\item The Avellaneda-Lin  type estimates were extended to the case of 
 Neumann boundary conditions by  Kenig,  Lin and  Shen  \cite{KLS1,KLS12,KLS3} (see also 
\cite{BLA12} for a related work). 
 These estimates should be helpful to study the next order approximation 
  for the  Neumann boundary condition case 
 
\item Many new probabilistic results were proved when an interface is present 
 (see \cite{HM11}) or in the trying to compute  the accurate value of the 
  homogenized matrix (see \cite{GO11,GO12}). 

\item Some different method was used to compute homogenized boundary 
data for none oscillating coefficient (\cite{LS,ASS}). 

 \end{itemize}

 

\def\cprime{$'$}


\begin{thebibliography}{10}




 \bibitem{ASS}  Aleksanyan, H. Shahgholian, H., Sj\"olin, P.:  Applications of Fourier analysis in homogenization and boundary layer. Available at arXiv:1205.5210v2 (2012)


\bibitem{AA99}
G.~Allaire and M.~Amar.
\newblock Boundary layer tails in periodic homogenization.
\newblock {\em ESAIM Control Optim. Calc. Var.}, 4:209--243 (electronic), 1999.





\bibitem{AvLi}
M.~Avellaneda and F.-H. Lin.
\newblock Compactness methods in the theory of homogenization.
\newblock {\em Comm. Pure Appl. Math.}, 40(6):803--847, 1987.

\bibitem{BLP78}
A.~Bensoussan, J.-L. Lions, and G.~Papanicolaou.
\newblock {\em Asymptotic analysis for periodic structures}.
\newblock North-Holland Publishing Co., Amsterdam, 1978.

\bibitem{BLA12}
X.~Blanc, F.~Legoll, and A.~Anantharaman.
\newblock Asymptotic behavior of green functions of divergence form operators
  with periodic coefficients.
\newblock {\em Applied Mathematics Research eXpress}, 2012.


\bibitem{CD99b}
D.~Cioranescu and P.~Donato.
\newblock {\em An introduction to homogenization}, volume~17 of {\em Oxford
  Lecture Series in Mathematics and its Applications}.
\newblock The Clarendon Press Oxford University Press, New York, 1999.

\bibitem{GM10}
D.~G{\'e}rard-Varet and N.~Masmoudi.
\newblock Relevance of the slip condition for fluid flows near an irregular
  boundary.
\newblock {\em Comm. Math. Phys.}, 295(1):99--137, 2010.

\bibitem{GM11}
D.~G{\'e}rard-Varet and N.~Masmoudi.
\newblock Homogenization in polygonal domains.
\newblock {\em J. Eur. Math. Soc. (JEMS)}, 13(5):1477--1503, 2011.

\bibitem{GM12}
D.~G{\'e}rard-Varet and N.~Masmoudi.
\newblock Homogenization and boundary layers.
\newblock {\em Acta Math.}, 209(1):133--178, 2012.



\bibitem{GO11}
A.~Gloria and F.~Otto.
\newblock An optimal variance estimate in stochastic homogenization of discrete
  elliptic equations.
\newblock {\em Ann. Probab.}, 39(3):779--856, 2011.

\bibitem{GO12}
A.~Gloria and F.~Otto.
\newblock An optimal error estimate in stochastic homogenization of discrete
  elliptic equations.
\newblock {\em Ann. Appl. Probab.}, 22(1):1--28, 2012.



\bibitem{HM11}
M.~Hairer and C.~Manson.
\newblock Periodic homogenization with an interface: the multi-dimensional
  case.
\newblock {\em The Annals of Probability}, 39(2):648--682, 2011.



\bibitem{JKO94}
V.~V. Jikov, S.~M. Kozlov, and O.~A. Oleinik.
\newblock {\em Homogenization of differential operators and integral
  functionals}.
\newblock Springer-Verlag, Berlin, 1994.
\newblock Translated from the Russian by G. A. Yosifian.



 \bibitem{KLS1} Kenig, C. , Lin, F., Shen, Zh.:  Periodic Homogenization of Green and Neumann Functions. Available at arXiv:1201.1440v1 (2012)

    \bibitem{KLS12}
C.~E. Kenig, F.~Lin, and Z.~Shen.
\newblock Convergence rates in {$L^2$} for elliptic homogenization problems.
\newblock {\em Arch. Ration. Mech. Anal.}, 203(3):1009--1036, 2012.

 
        \bibitem{KLS3} Kenig, C. , Lin, F., Shen, Zh.: Homogenization of Elliptic Systems with Neumann Boundary Conditions. Available at 	arXiv:1010.6114v1 (2010)

        \bibitem{LS} Lee, K., Shahgholian, H.: Homogenization of the boundary value for the Dirichlet problem. Avaliable at arXiv:1201.6683v1 (2012)



\bibitem{JLL}
J.-L. Lions.
\newblock {\em Some methods in the mathematical analysis of systems and their
  control}.
\newblock Kexue Chubanshe (Science Press), Beijing, 1981.

\bibitem{MV97}
S.~Moskow and M.~Vogelius.
\newblock First-order corrections to the homogenised eigenvalues of a periodic
  composite medium. {A} convergence proof.
\newblock {\em Proc. Roy. Soc. Edinburgh Sect. A}, 127(6):1263--1299, 1997.

\bibitem{MT97}
F.~Murat and L.~Tartar.
\newblock Calculus of variations and homogenization [ {MR}0844873 (87i:73059)].
\newblock In {\em Topics in the mathematical modelling of composite materials},
  volume~31 of {\em Progr. Nonlinear Differential Equations Appl.}, pages
  139--173. Birkh\"auser Boston, Boston, MA, 1997.


\bibitem{Prange1}
C. Prange 
\newblock First-order expansion for the Dirichlet eigenvalues of an elliptic system with oscillating coefficients
\newblock A paraître dans {\em Asymptotic Analysis}

\bibitem{Prange2}
C. Prange 
\newblock Asymptotic analysis of boundary layer correctors in periodic homogenization 
\newblock A paraître dans {\em SIAM Journal on Mathematical Analysis}





\bibitem{Sanchez80}
E.~S{\'a}nchez-Palencia.
\newblock {\em Nonhomogeneous media and vibration theory}.
\newblock Springer-Verlag, Berlin, 1980.













\end{thebibliography}
 


\end{document}